\documentclass[11pt]{amsart}

\usepackage{amsthm,amssymb,amsmath,amsopn,amsfonts,etoolbox,mathptmx}
\usepackage{stmaryrd,calc,mathrsfs,enumerate}
\usepackage{tikz}\usetikzlibrary{matrix,arrows} 

\makeatletter
\patchcmd{\@settitle}{\uppercasenonmath\@title}{}{}{}
\patchcmd{\@setauthors}{\MakeUppercase}{\scshape}{}{}
\patchcmd{\section}{\scshape}{\bfseries}{}{}
\renewcommand{\@secnumfont}{\bfseries}
\patchcmd{\abstract}{\scshape\abstractname}{\textbf{\abstractname}}{}{}
\makeatother

\usepackage{hyperref} 
\hypersetup{
	colorlinks=true,
	citecolor=blue,
	linkcolor=red
}

\newtheorem{thm}{Theorem}[section]
\newtheorem{prop}[thm]{Proposition}
\newtheorem{lem}[thm]{Lemma}
\newtheorem{cor}[thm]{Corollary}

\theoremstyle{definition}
\newtheorem{definition}[thm]{Definition}
\newtheorem{example}[thm]{Example}

\theoremstyle{remark}
\newtheorem{remark}[thm]{Remark}

\numberwithin{equation}{section}

\newcommand{\bdef}{\begin{definition}}\newcommand{\ndef}{\end{definition}}
\newcommand{\bteo}{\begin{thm}}\newcommand{\nteo}{\end{thm}}
\newcommand{\bprop}{\begin{prop}}\newcommand{\nprop}{\end{prop}}
\newcommand{\brmk}{\begin{remark}}\newcommand{\nrmk}{\end{remark}}
\newcommand{\bcor}{\begin{cor}}\newcommand{\ncor}{\end{cor}}
\newcommand{\blem}{\begin{lem}}\newcommand{\nlem}{\end{lem}}
\newcommand{\bex}{\begin{example}}\newcommand{\nex}{\end{example}}
 \newcommand{\Z}{\mathbb{Z}}

\newcommand{\Neron}{N\'{e}ron} 

\renewcommand{\P}{\mathbb{P}}
\renewcommand{\O}{\mathcal{O}}

\newcommand{\g}{\mathcal}

\newcommand{\f}{\mathfrak}
\newcommand{\qt}[1]{``#1''}

\newcommand{\nik}{\mathcal{F}_g^{\mathfrak{N}}}

\newcommand*{\dashedarrow}[1][]{\mathbin{\tikz [baseline=-0.25ex,-latex, dashed,#1] \draw [#1] (0pt,0.5ex) -- (1.5em,0.5ex);}}

\DeclareMathOperator{\rk}{rk}
\DeclareMathOperator{\Pic}{Pic}

\DeclareMathOperator{\NS}{NS}

\DeclareMathOperator{\Sing}{Sing}

\DeclareMathOperator{\Cliff}{Cliff}

\DeclareMathOperator{\Aut}{Aut}

\DeclareMathOperator{\codim}{codim}

\DeclareMathOperator{\td}{td}
\DeclareMathOperator{\ch}{ch}
\DeclareMathOperator{\CH}{CH}

\begin{document}


\title{A note on Nikulin surfaces and their moduli spaces}
\author{Marco Ramponi}
\address{Lehrstuhl f\"{u}r Algebra und Zahlentheorie, Institut f\"{u}r Mathematik - Universit\"{a}t Augsburg, D-86135 Augsburg, Deutschland}


\begin{abstract}
We study a number of  natural linear systems carried by any polarized Nikulin surface of genus $g$, determine their positivity and establish their Brill-Noether theory. As an application, we compute the classes of some natural effective divisors on the moduli space of Nikulin surfaces, relying upon recent work of Farkas and Rim\'{a}nyi.


\end{abstract}

\maketitle

\section{Introduction}

A polarized Nikulin surface of genus $g\geq2$ is a smooth polarized K3 surface $(X,L)$, where $L$ is a big and nef line bundle with self-intersection $L^2=2g-2$, such that there exists a set of disjoint smooth rational curves $R_1,\dots,R_{8}$ on $X$ satisfying $L\cdot R_i =0$ for each $i=1,\dots,8$. 
The divisor class of $R:= R_1+\dots+R_{8}$ is divisible by $2$ in the Picard group of $X$ and, denoting by $e$ the primitive class satisfying 
\[ e^{\otimes 2} = \O_X(R) , \] 
one can describe the Nikulin lattice $\f{N}$ as the lattice generated by the classes of $R_1,\ldots,R_8$ and $e$. This yields a primitive embedding 
\[\Lambda_g := \Z\cdot [L] \oplus \f{N} \hookrightarrow \Pic(X) , \] 
from which one can construct the moduli space $\nik$ of Nikulin surfaces of genus $g$, which is irreducible and $11$-dimensional \cite{Dol96,vGS07}.

Nikulin surfaces represent a rather special class of K3 surfaces, which have been studied in relation with various topics, including the theory of automorphisms \cite{Nik76,GS08}, moduli spaces \cite{Mor84,vGS07}, the study of Prym curves \cite{FK16} and of the birational geometry of their moduli spaces \cite{FV12,FV16,KLCV,Ver16}. 

Following \cite{GS08}, the very general point of $\nik$ corresponds to a polarized Nikulin surface $(X,L)$ with Picard number nine and, if $g\geq3$, the line bundle $L$ induces a morphism which contracts the $8$ curves $R_i$ and maps $X$ to a surface in $\P^g$ with $8$ ordinary double points. If $g=4$ or $5$, the line bundle $L\otimes e^{-1}$ is ample, and very ample for $g\geq6$, mapping $X$ to a surface in $\P^{g-2}$. These linear systems and their interplay provide a source of interesting geometry which has been proven useful for several effective constructions in some of the aforementioned works. A detailed description of these projective models has been carried out for Nikulin surfaces of low genus in \cite{GS08}.

The starting point of this work is a systematic approach for the study of these linear systems on Nikulin surfaces of any genus. Consider the series of line bundles
\[ L_m = L \otimes e^{-m}, \quad m=0,1,2,\ldots  \]

Assuming $g$ is large enough with respect to $m$, how positive is $L_m$? For very general $(X,L)$ in $\nik$, the situation is as good as one could expect. 

\bteo\label{mapsLm} Let $(X,L)$ be a Nikulin surface of genus $g$ with Picard number nine. Write $g=2k^2+ p $, where $k\geq1$ and $0\leq p  < 4k+2$. Denote by 
\[ L_m = L \otimes e^{-m} , \quad m=0,1,\ldots,k. \]
Then,
\begin{enumerate}[(i)]
\item For any $1\leq m\leq k-1$ (in particular $g\geq 8$) the general member of the linear system $|L_m|$ is a smooth irreducible curve of genus 
\[g_m = g - 2m^2 \geq 6. \] 
In fact, $L_m$ is very ample and defines an embedding of $X$ in $\P^{g_m}$. 
\item $|L_k|$ contains a smooth irreducible curve of genus $ p $. Moreover, $L_k$ is ample for $ p =2$, and very ample for $ p \geq3$.
\end{enumerate}
\nteo

For example, one easily checks that for any $1\leq m\leq k-1$, a smooth irreducible curve $D$ in the linear system $|L_{m+1}|$ is embedded (non-specially) by $L_m$ as a non-degenerate linearly normal curve of degree $2(g_m -2m-1)$ in $\P^{g_m}$. In particular, letting $g=2k^2+ p $,  the curves $D\in|L_k|$ are embedded by $L_{k-1}$ as linearly normal curves of genus $ p $ in $\P^{4k-2+ p }$. So, for example:
\begin{itemize}
\item $ p =0$ ($g=2k^2$) : $D$ is the rational normal curve in $\P^{4k-2}$,
\item $ p =1$ ($g=2k^2+1$) : $D$ is a normal elliptic curve in $\P^{4k-1}$,
\item $ p =2$ ($g=2k^2+2$) : $D$ is a normal genus $2$ curve in $\P^{4k}$,
\end{itemize}
and so on. Let us point out that in genus $g=8$ (i.e.\,the case $k=2$ and $ p =0$), the existence of the smooth rational curve $D \in|L_2|$ embedded by $L_1$ as a rational normal curve of degree $6$ in $\P^6$ is the starting point of Verra's construction\,\cite{Ver16} for the proof of the rationality of $\g{F}_8^{\f{N}}$.

Returning to the linear systems $|L_m|$, the next natural step is to determine the Brill-Noether theory of the hyperplane sections of $X$ in these projective models. 

\bteo\label{BNgeneral} Under the same assumptions and notations of Theorem \ref{mapsLm}, all smooth curves in $|L_m|$, for any $0\leq m\leq k$, are Brill-Noether general, i.e.\,have maximal Clifford index $\lfloor \frac{g_m-1}{2} \rfloor$.
\nteo

We note that for $m=1$ the statement of Theorem \ref{BNgeneral} was proved by Farkas and Kemeny in their proof of the Prym-Green conjecture in odd genus \cite{FK16}.

In the second part of the paper we use the projective geometry of Nikulin surfaces in order to compute the classes of some natural effective divisors in $\nik$. This was directly inspired by the recent paper of Farkas and Rim\'{a}nyi \cite{FR}, and relies on some beautiful general formulas therein contained. Specifically, for any index $m=0,\ldots,k$ (assuming $p\geq 3$ when $m=k$), we consider the virtual divisors

\[  \f{D}^{\rk (4)}_m = \{ (X,L) \in \nik \, \colon \, \, I_m^{\rk(4)} (2) \neq 0  \} , \]

where $I_m(2)$ is the vector space of quadrics containing the image of $X$ by the morphism $X \to \P^{g_m}$ induced by $L_m$ and 
\[ I_m^{\rk (4)}(2) = \{ q \in I_m(2) \colon \, \rk(q) \leq 4 \} . \] 

Along the lines of \cite[\S 4]{MOP17}, we consider the universal Nikulin surface
\[ \pi \colon \g{X} \longrightarrow \nik \]
together with a choice of universal bundles $\g{L}$ and $\g{E}$ associated to the polarization $L$ and the class $e$ of each fiber $(X,L)$. We define the tautological classes

\[ \kappa_{a,b,c} = \pi_\ast \left( c_1(\g{L})^a \cdot c_1(\g{E})^b \cdot c_2( \g{T}_\pi )^c \right) \in \mathrm{CH}^{a+b +2c -2} (\nik) , \]
for some non-negative integers $a,b,c$. In codimension one, the six classes
\[ \kappa_{3,0,0}, \, \, \kappa_{0,3,0}, \, \, \kappa_{2,1,0}, \, \, \kappa_{1,2,0}, \, \, \kappa_{1,0,1}, \, \, \kappa_{0,1,1}  \]
depend on the choice of $\g{L}$ and $\g{E}$. However, the following linear combinations are independent from this choice, and thus yield well defined classes in $\mathrm{CH}^1(\nik)$:

\begin{align*}
\gamma_0 = \kappa_{3,0,0} - \frac{(g-1)}{4} \kappa_{1,0,1} \, , \quad \gamma_2 = \kappa_{1,0,1} + 6 \kappa_{1,2,0} \, , \\
\gamma_1 = \kappa_{2,1,0} - \frac{(g-1)}{12} \kappa_{0,1,1} \, , \quad \gamma_3 =\kappa_{0,1,1} + 2 \kappa_{0,3,0} \, .
\end{align*}
(cf.\,\S \ref{tauto} below for more details). It turns out that the class of $\f{D}^{\rk (4)}_m$ can be computed in terms of these four classes $\gamma_0, \gamma_1, \gamma_2, \gamma_3 \in \CH^1(\nik)$, together with the Hodge class $\lambda \in\CH^1(\nik)$. Precisely, we have:

\bteo\label{THM-class-of-Dm}
Let $g=2k^2+ p $, for some integers $k\geq1$ and $0\leq p  < 4k+2$. If $m=0,\ldots,k-1$ or $m=k$ and $p\geq 3$, then $\f{D}^{\rk (4)}_m$ is an effective divisor in $\nik$ and its class in the Chow ring $\CH^1(\nik)$ can be expressed as
\[ \left[ \f{D}^{\rk (4)}_m \right] =  A_m \left( \frac{ 2 }{ g_m +1 } \gamma_0 - \frac{ 6m }{ g_m +1 } \gamma_1 + \frac{ m^2 }{ g_m +1 } \gamma_2 - \frac{ m^3 }{ g_m +1 } \gamma_3 + (2g_m-1) \lambda \right) , \]
where $g_m = g-2m^2$ and some rational number $A_m$. 
\nteo

The proof rests upon a Grothendieck-Riemann-Roch calculation on the universal Nikulin surface of genus $g$ and the application of general formulas from \cite{FR} which allow one to express $[\f{D}^{\rk (4)}_m]$ in terms of Chern classes of vector bundles on $\nik$.

\section*{Acknowledgements}
This work started thanks to the financial support of the Einstein Foundation and I wish to thank G.\,Farkas and R.\,Pandharipande for giving me the opportunity to spend some pleasant time in Berlin. I am warmly thankful to G.\,Farkas for useful conversations around these topics.

\section*{Notation and conventions}

Nikulin surfaces of genus $g$, as defined in the Introduction, exist in all genera and are sometimes called \emph{standard} Nikulin surfaces. In odd genus (only) there exist also Nikulin surfaces for which $L\cdot e\neq0$, i.e.\,one does not have an orthogonal sum decomposition $\Z\cdot [L] \oplus \f{N}$. Accordingly, the \Neron-Severi lattice of the very general such K3 surface is rather an \emph{overlattice of index $2$} of $\Lambda_g$, cf.\,\cite{vGS07}. In this paper we are only concerned with standard Nikulin surfaces, whence we will omit this specification throughout. We work over the field of complex numbers.

\section{Polarized Nikulin surfaces}

We briefly recall some basic facts on Nikulin surfaces and refer to \cite{GS08} for details. Let $Y$ be a K3 surface carrying a symplectic involution $\iota\in\Aut(Y)$. Then $\iota$ has $8$ isolated fixed points and the quotient $Y/\iota$ is a K3 surface with eight ordinary double points. Let $\sigma\colon \widetilde{Y} \to Y$ be the blow-up of $Y$ at the eight fixed points of $\iota$. The involution $\iota$ naturally lifts to an involution on $\widetilde{Y}$, fixing the eight exceptional divisors $E_1,\ldots,E_8$. We denote the quotient surface by this involution by $X$ (which turns out to be a K3 surface) and we let $f\colon \widetilde{Y} \to X$ be the quotient map. 

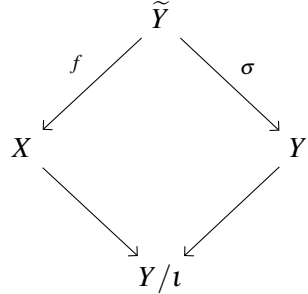
\[\begin{tikzpicture}[>=angle 90]
\matrix(a)[matrix of math nodes,
row sep=3em, column sep=3em,
text height=1.5ex, text depth=0.25ex]
{ & \widetilde{Y} & \\
X & & Y \\
 & Y / \iota & \\};
\path[->,font=\scriptsize]
(a-1-2) edge node[above left]{$f$} (a-2-1)
(a-1-2) edge node[above right]{$\sigma$} (a-2-3)
(a-2-1) edge (a-3-2)
(a-2-3) edge (a-3-2);
\end{tikzpicture}\] 

By construction, the ramification divisor of $f$ is the sum of the rational curves $R_i := f(E_i)$, and is therefore $2$-divisible in $\Pic(X)$. We denote by $e$ the primitive class $e\in\Pic(X)$ satisfying the identity
\[ e^{\otimes 2} = \O_X(R) , \]
where $R= R_1+\dots+R_{8}$. The classes of $R_1,\dots,R_{8}$ and $e$ span the Nikulin lattice $\f{N}$, and one gets in this way a primitive embedding $\f{N} \hookrightarrow \Pic(X)$.

\bdef
A \emph{ polarized Nikulin surface of genus $g\geq2$ } consists of a K3 surface $X$ and a primitive embedding $j\colon \Lambda_g \hookrightarrow\Pic(X)$, where $\Lambda_g := \Z\cdot L \oplus \f{N}$, such that $L^2=2g-2$ and $j(L)$ is a big and nef class. With a small abuse of notation, we denote a polarized Nikulin surfaces simply by the pair $(X,L)$.
\ndef

Polarized Nikulin surfaces of genus $g$ form a moduli space denoted by $\nik$, which is  known to be irreducible and $11$-dimensional \cite{Dol96,vGS07}. The very general point of $\nik$ corresponds to a K3 surface $X$ with Picard number nine. More precisely, one has an isomorphism $\Pic(X)\simeq\Lambda_g$, see \cite[Prop.\,2.1]{GS08}.

For the rest of this section we assume $(X,L)$ to be a polarized Nikulin surface of genus $g$ with Picard number nine and write 
\[ g=2k^2+ p , \] 
where $k\geq1$ and $0\leq p  < 4k+2$. We denote by
\[ L_m = L \otimes e^{-m} .\]

\blem $L_m$ is ample if $1\leq m\leq k-1$ or $m=k$ and $ p \geq2$.
\nlem

\begin{proof}
Since $L_m$ is effective and $(L_m)^2>0$, we have that $L_m$ is ample if and only if $L_m$ intersects positively any smooth rational curve on $X$. Let $D$ be an effective divisor and write $D \sim aL -b_1 R_1-\cdots -b_8R_8$, where $a\in\Z$ and $b_i\in\frac{1}{2}\Z$. Then,
\[ D\cdot L_m = 2a(g-1) - m(\sum_{i=1}^8 b_i). \]
Assuming $D\cdot L_m\leq0$, we are going to show that $D^2<-2$ for $1\leq m\leq k-1$. Since $L_m$ intersects positively each $R_i$, we can assume that $D$ does not contain $R_i$ as a component, whence $D\cdot R_i=2b_i\geq0$. Hence we must have $a>0$, else $D$ would not be effective. Thus $D\cdot L_m\leq0$ yields $2a(g-1)\leq m(b_1 +\ldots+ b_8)$, and we can square both sides of this inequality, obtaining
\[ 4a^2(g-1)^2 \leq m^2(\sum_{i=1}^8 b_i)^2 \leq 8m^2 \sum_{i=1}^8 b_i^2. \]
where the latter estimate follows by the Cauchy-Schwarz inequality. Let us re-write this condition as
\begin{align} \label{cauchy-schwarz1}
\frac{a^2(g-1)^2}{m^2} \leq 2 \sum_{i=1}^8 b_i^2. 
\end{align}

Assuming the latter inequality, we can estimate the self-intersection of $D$,
\begin{align*}
D^2 & = 2a^2(g-1) - 2\sum_{i=1}^8 b_i^2 \\
& \leq 2a^2(g-1) - \frac{a^2(g-1)^2}{m^2}  \quad \mbox{ by (\ref{cauchy-schwarz1}) } \\
& = a^2 \left(\frac{ g-1 }{ m^2 } \right) (1-g_m).
\end{align*}

Since $1-g_m <0$, in order for the condition $D^2 < -2$ to be satisfied, it is enough to ask for the following numerical condition:
\[  a^2 \left(\frac{ g-1 }{ m^2 } \right)(g_m-1) > 2. \]
By assumption $g_m=g-2m^2\geq 6$, or equivalently $m^2 \leq \frac{1}{2}(g-6)$, and using both these estimates we finally get
\[  a^2 \left(\frac{ g-1 }{ m^2 } \right)(g_m-1) \geq a^2 \cdot 2 \left(\frac{ g-1 }{ g-6 } \right) \cdot 5 > 10 a^2 > 2. \]
This proves that $L_m$ is ample for $1\leq m\leq k-1$. Moreover, we observe that the same argument works for $m=k$ and $ p \geq2$, in which case $g_k =  p $ and $k^2=\frac{1}{2}(g- p )$, so that we still get
\[  a^2 \left(\frac{ g-1 }{ k^2 } \right)(g_k-1) = a^2  \cdot 2 \left(\frac{ g-1 }{ g- p  } \right) ( p -1) > 2. \]

This concludes the proof. 
\end{proof}

\blem $L_m$ is very ample for any $1\leq m\leq k-1$ or $m=k$ and $ p \geq3$.
\nlem

\begin{proof}
So far we have determined the amplitude of $L_m$ for any $1\leq m\leq k$, with $ p \geq2$ for $m=k$. Since each such $L_m$ is big and nef, by the results of Saint-Donat \cite{SD} there can be only two types of obstruction for $L_m$ to be very ample: either the existence of fixed components, or of hyperelliptic curves in $|L_m|$. Specifically, we have to exclude the following three possibilities:
\begin{enumerate}[(a)]
\item There exist an elliptic curve $E$ and a smooth rational curve $\Gamma$ on $X$ such that $L_m=\O_X(rE+\Gamma)$, with $r\geq2$ and $E\cdot\Gamma=1$. 
\item There exists an elliptic curve $E$ such that $L_m\cdot E =2$. 
\item There exists a curve $B$ of genus $2$ such that $L_m=\O_X(2B)$.
\end{enumerate}

Notice that (c) may occur in genus $g_m=5$ only, but in fact it is automatically excluded by the properties of the \Neron-Severi lattice of the very general Nikulin surface $X$, since the class $L_m=L\otimes e^{-m}$ is not $2$-divisible in $\NS(X)$ \cite[Prop.\,2.1]{GS08}. 

We will proceed by excluding cases (a) and (b) simultaneously, by showing that for any elliptic curve $E$ we have $L_m \cdot E >2$. Let us sketch the argument: given an effective divisor $D \sim aL -b_1 R_1-\cdots -b_8R_8$, we assume $D\cdot L_m\leq 2$. This yields the inequality
\begin{align*}
\frac{ [a(g-1)-1]^2 }{m^2} \leq 2 \sum_{i=1}^8 b_i^2. 
\end{align*}

Assuming the latter inequality, we can estimate the self-intersection of $D$,
\begin{align*}
D^2 & = 2a^2(g-1) - 2\sum_{i=1}^8 b_i^2 \\
& \leq \frac{1}{m^2} [ a^2 (g-1) (1-g_m) + 2a(g-1)-1 ].
\end{align*}
We use the latter expression to impose the condition $D^2 > 0$. This yields to a quadratic equation in $a$, which will always be satisfied for any

\[ a > \frac{1}{g-2m^2-1} \left( \sqrt{ \frac{g-1+2m^2}{2g-2} } + 1 \right).  \]

Since $a\geq 1$, the latter inequality will hold for any $a$ as long as we impose the right hand side to be strictly lower than 1. Hence, we study the following condition:

\[ f(T):=\frac{1}{g-T-1} \left( \sqrt{ \frac{g-1+T}{2g-2} } + 1 \right) < 1 \quad  (T=2m^2). \]

By direct computation, we find $f(T)<1$ for any $g>2$ and
\[ 0 < T < \frac{2g^2-8g+7}{2g-2} .\]

Finally, since $T=2m^2 \leq 2k^2 = g- p $, we are led to the inequality
\[ g- p  < \frac{2g^2-8g+7}{2g-2} ,\]
which holds for all $g>3$ and $ p \geq3$ (and fails for $ p =2$, not surprisingly). This shows that $L_m$ is very ample and concludes the proof.
\end{proof}

Let once again $g=2k^2+ p$, where $k\geq1$ and $0\leq p  < 4k+2$. We now focus our attention on the line bundle $L_k = L \otimes e^{-k}$.

\blem $|L_k|$ contains a smooth irreducible curve of genus $ p \geq0$.
\nlem

\begin{proof}
By Riemann-Roch, $h^0(X,L_k)\geq 1+ p \geq 1$, whence $L_k$ is effective. Let
\[ |L_k| = |M| + F \]
be the moving and fixed part decomposition of the linear system $|L_k|$. We show that either $F=0$ (whence $|L_k|$ is basepoint free and the statement follows by Bertini's Theorem) or $M=0$ and $F$ is a smooth rational curve. 

When $k=1$ the statement follows by Proposition 3.2 and Proposition 3.4 in \cite{GS08}. Therefore, we assume $k\geq 2$. Then $L_{k-2}$ is effective and $h^1(X,L_{k-2})=0$. Consider the short exact sequence
\[ 0 \to L_k \to L_{k-2} \to \O_R(L_{k-2}) \to 0, \]
where $\O_X(R)=e^{\otimes 2}$. Note that $L_{k-2} \cdot R_i>0$ and the induced restriction map 
\[ H^0(X,L_{k-2}) \to H^0(R,\O_R(L_{k-2})) = \bigoplus_{i=1}^8 H^0(R_i,\O_{R_i}(L_{k-2})) \] 
is surjective. Then $h^1(X,L_k)=0$ and $h^0(X,L_k) = 1+ p $, by Riemann-Roch. 

If $ p \geq1$, the linear system $|L_k|$ has then a non-trivial moving part, and the established vanishing $h^1(X,L_k)=0$ implies that the general member of $|L_k|$ is irreducible, i.e. $F=0$ and $|L_k|$ is basepoint free. 
If $ p =0$, then $M=0$ and the linear system $|L_k|$ consists of a single effective divisor $F$ of self-intersection $F^2=-2$. If $F$ contains one of the $R_i$'s as a component, then $F\cdot R_i<0$. On the other hand, 
\[ F\cdot R_i = L_k\cdot R_i = -\frac{k}{2} (\sum R_j) \cdot R_i = k > 0 . \]
Thus, we write $F = \Gamma_1+\dots+\Gamma_r$ as a sum of its irreducible components, which are smooth rational curves each of which is distinct from the $R_i$'s. By the structure of $\Pic(X)$, we can write
\[ \Gamma_\ell = a_\ell \cdot L - \sum_{i=1}^8 b_{\ell,i} \, R_i. \]
For each $1\leq \ell \leq r$ and $1\leq i\leq 8$, we have $b_{\ell,i} \geq 0$, since $\Gamma_\ell$ is different from $R_i$. Therefore $a_\ell >0$, else $\Gamma_\ell$ would not be effective, and by $F=\Gamma_1+\dots+\Gamma_r$ and $F \sim L -k e$, we get $a_1+\dots+a_r = 1$, whence $r=1$.
\end{proof}

\begin{proof}[Proof of Theorem \ref{mapsLm}]
It follows at once by the previous lemmas.
\end{proof}

In the next section we focus our attention on the Brill-Noether aspects of the linear systems $|L_m|$ and show how this naturally leads to divisors in $\nik$.

\section{Brill-Noether general curves and the locus $\f{D}^{\rk (4)}_m$ }\label{sec-BN+quadric-loci}

Let us briefly recall some well-known facts from Brill-Noether theory. Let $C$ be a smooth algebraic curve of genus $g$. For any $A \in \Pic(C)$, we let
\[ \Cliff(A) = \deg A-2h^0(A)+2 . \] 
The Clifford index of $C$ is by definition 
\[ \Cliff(C) = \min\{\Cliff(A) \, \colon \, A\in\Pic(C), \, h^0(A)\geq2, \, h^1(A)\geq2\}. \]

Line bundles $A$ satisfying the conditions as in the definition of $\Cliff(C)$ are said to \emph{contribute to the Clifford index} of $C$. Such line bundles exist for any curve as long as $g\geq4$. When $g<4$ we adopt the standard convention 
\[ \Cliff(C) := \begin{cases} 0 \quad \mbox{ if } \, g=1,2 \, \mbox{ or 3 and $C$ is hyperelliptic} \\  1 \quad \mbox{ if } \, g=3 \, \mbox{ and $C$ is non-hyperelliptic} \end{cases} \]

For any curve $C$ of genus $g$, one has the inequality
\[ \Cliff(C)\leq \left\lfloor \frac{g-1}{2} \right\rfloor , \]
which is an equality for the general curve  of genus $g$. 

\bdef Let $S$ be a surface. We say that a curve $C\subset S$, or equivalently the line bundle $\O_S(C)$, decomposes into a sum of movable classes if $C$ is linearly equivalent to a sum $D_1+D_2$ of two divisors satisfying $h^0(D_i) \geq 2$, for $i=1,2$.
\ndef

For the rest of this section, we let $(X,L) \in \nik$ be a polarized Nikulin surface of genus $g$. We write $g=2k^2+ p $, where $k\geq1$ and $0\leq p < 4k+2$ and denote by $L_m = L \otimes e^{-m}$, keeping the same notations of the previous section.

\bprop\label{No-decompositions} Assume $X$ has Picard number $9$. If $m=0,\ldots,k-1$ or $m=k$ and $p\geq 2$, the line bundle $L_m$ does not decompose into a sum of movable classes.
\nprop

\begin{proof}
Let us assume by contradiction that there is some smooth curve $D\in|L_m|$ which decomposes into a sum of movable classes. We can clearly assume $D^2 \geq2$. Among all such decompositions, we choose one $D\sim D_1+D_2$ such that the intersection $D_1\cdot D_2$ is minimal. We can then assume that one of the two classes, say $D_1$, is base point free and the general member of $|D_1|$ is a smooth irreducible curve, cf.\,\cite[Prop\,2.7]{Kn04}. In particular, $D_1$ is a nef class, and the restriction $(D_1)|_D$ contributes to the Clifford index of $D$, whence $\deg_D(D_1)>0$. We want to show that $D_2$ is then not a movable class. Up to linear equivalence, we can write
\[ D_\ell = a_\ell \cdot L - \sum_{i=1}^8 b_{\ell,i} \, R_i, \quad \ell=1,2, \]
where $a_\ell\in\Z$ and $2 b_{\ell,i} \in \Z$.
By intersecting $D_\ell$ with the nef class $L$, we get $a_\ell\geq0$ for $\ell=1,2$. Since $D_1$ is nef, $D_1\cdot R_i = 2b_{1,i} \geq0$ for all $i=1,\dots,8$. Finally,
\[ D \cdot D_1 = (L -me) \cdot D_1= a_1 L^2 - m\sum_{i=1}^8 b_{1,i} > 0 , \]
which implies $a_1>0$. On the other hand, by $D \sim D_1+D_2$ we have $a_1+a_2=1$. Therefore $a_2=0$, which clearly contradicts $h^0(D_2)\geq 2$.
\end{proof}

We can now prove Theorem \ref{BNgeneral} from the Introduction.

\begin{proof}[Proof of Theorem.\,\ref{BNgeneral}]
Let $D\in|L_m|$ be a smooth curve. Assuming by contradiction that $D$ has Clifford index $\Cliff(D)< \lfloor\frac{g_m-1}{2} \rfloor$,  there exists a line bundle $M$ on $X$ such that $M\otimes\O_D$ contributes to the Clifford index of $D$ and $\Cliff(D)=\Cliff(M\otimes\O_D)$, by \cite{GL}.
By the definition of Clifford index it follows that $h^0(M)\geq2$ and also $h^0(D-M)\geq2$. Thus $D \sim D_1 + D_2$, with $D_1\in |M|$ and $D_2\in |D -M|$ is a decomposition of $D$ into two movable classes, contradicting Proposition \ref{No-decompositions}.
\end{proof}

We are now going to define the divisors $\f{D}^{\rk (4)}_m$ from the Introduction. The following discussion parallels \cite[\S 9.1]{FR}.
Consider an element $(X,L)$ in $\nik$ and for any $0 \leq m \leq k$ we denote by $L_m = L \otimes e^{-k}$. Each one of the following conditions,
\begin{enumerate}[(a)]
\item There exists an elliptic pencil $E$ on $X$ such that $E\cdot L_m = 1$,
\item There exists an elliptic pencil $E$ on $X$ such that $E\cdot L_m = 2$,
\item There exists a smooth rational curve $R$ on $X$ such that $R\cdot L_m = 0$,
\end{enumerate}
singles out a virtual Noether-Lefschetz divisor in $\nik$, i.e.\,can be rephrased as the condition of the existence of a primitive embedding $\Lambda_g \hookrightarrow \Lambda_g'$, for some lattices $\Lambda_g'$ of rank $10$. \qt{Virtual} here means that these conditions may possibly be empty on $\nik$. However, up to assuming $p\geq3$ in the case $m=k$, Theorem \ref{mapsLm} and Theorem \ref{BNgeneral} guarantee that for general $(X,L)$ the conditions (a) and (b) are not satisfied, whence they define an actual divisor in $\nik$. Moreover, condition (c) defines an actual divisor as long as $m\geq1$.

Well-known results of Saint-Donat imply that Nikulin surfaces $(X,L)$ outside the union of these divisors are such that each $L_m$ is basepoint free, (pseudo-) ample and non-hyperelliptic, whence the multiplication map
\[ \psi_m \colon S^2 H^0(X, L_m) \to H^0(X, L_m^{\otimes 2} ) \]
is surjective, cf.\,\cite[Thm.\,6.1]{SD}. Thus $I_m(2) := \ker(\psi_m)$ in $S^2 H^0(X,L_m)$ has codimension $h^0(X,L_m^{\otimes 2}) = 4g_m -2$, by Riemann-Roch. The closed subscheme

\[ \Sigma^{\rk(k)}_m = \{ q \colon \, \rk(q) \leq k\} \subset S^2 H^0(X, L_m) \]
has codimension 

\[ \binom{g_m+2-k}{2} , \] 
see e.g.\,\cite{Ott95}. Therefore, the expected codimension of scheme-theoretic intersection

\begin{align*}
I_m^{\rk (k)}(2) := \, \, I_m(2) \cap \Sigma_m^{\rk(k)} = \, \{ q \in I_m(2) \colon \, \rk(q) \leq k \} 
\end{align*}
is
\[ \codim I_m(2) + \codim \Sigma_m^{\rk(k)} = 4g_m - 2  + \binom{g_m+2 -k}{2} \]
in $S^2 H^0(X, L_m)$. In other words, we expect

\[ \dim I_m^{\rk(k)}(2) = (k-4) g_m + \frac{1}{2} (4+3k -k^2) ,\]
which vanishes precisely when $k=4$. We thus expect that for the general Nikulin surface $ I_m^{\rk(k)}(2) = 0$, so that, interpreting $\psi_m$ as a morphism between vector bundles over $\nik$, it follows that the condition $I_m^{\rk(4)} (2) \neq 0$ is a divisorial --or eventually empty-- condition on $\nik$, i.e.\,the locus

\[  \f{D}^{\rk (4)}_m = \{ (X,L) \in \nik \, \colon \, \, I_m^{\rk(4)} (2) \neq 0  \} , \]
is a virtual divisor in $\nik$. 

\bprop Let $0\leq m\leq k$ and assume $p\geq3$ if $m=k$. The locus $ \f{D}^{\rk (4)}_m$ is a Noether--Lefschetz divisor on $\nik$ supported on the subset of $\nik$ defined by the condition that $L_m$ decomposes into a sum of two movable classes.
\nprop
\begin{proof}
First we observe that $ \f{D}^{\rk (4)}_m \neq \nik$, by Proposition \ref{No-decompositions}. Assume $(X,L)$ is such that $L_m$ decomposes as a sum of movable classes $D_1,D_2$. Choose a general 2-dimensional subspace $V_i$ of $H^0(X, D_i)$ for $i=1,2$. The multiplication map $V_1 \otimes V_2 \rightarrow H^0(X, L_m)$ yields a $2 \times 2$ matrix whose determinant vanishes on the image of $X$ in $\P H^0(X, L_m)$, thus producing a rank four quadric containing $X$.
Conversely, suppose $X\subset\P^{g_m}$ is contained in a quadric $Q$ of rank at most $4$. We can assume $\rk Q = 4$, whence $\Sing Q = \P^{g_m-4} \subset\P^{g_m}$ and $Q$ is the cone over $\P^1\times\P^1$ under the projection $\P^{g_m} \dashedarrow \P^3$ from $\Sing Q$. The two rulings of $\P^1\times \P^1$ pull back to $Q$ and cut out on $X$ two divisors $D_1, D_2$, with $h^0(X,D_i)\geq2$, such that $L \sim D_1 + D_2$.
\end{proof}

\section{Computation of the class of $\f{D}^{\rk (4)}_m$ }\label{tauto}

Let $\nik$ be the moduli stack of polarized Nikulin surfaces and
\[ \pi \colon \g{X} \longrightarrow \nik \]
the universal polarized Nikulin surface of genus $g$. Denote by 
\[ \g{L}\in\Pic(\g{X}) \] 
the universal polarization. The Hodge bundle $\mathbb{E}$ on $\nik$ is defined by 
\[ \mathbb{E}^\vee = R^2 \pi_\ast \O_\g{X} \]
and we denote its first Chern class by 
\[ \lambda = c_1(\mathbb{E}) \in \mathrm{CH}^1(\nik) ,\]
which we refer to as the Hodge class on $\nik$. 

The relative tangent bundle $\g{T}_\pi \longrightarrow \g{X}$ is defined by the short exact sequence
\[ 0 \to \g{T}_\pi \to \g{T}_\g{X} \to \pi^\ast \g{T}_{\nik} \to 0 . \]
The relative canonical bundle $\omega_\pi = \det(\g{T}_\pi)^{-1}$, which restricts to the canonical (whence trivial) line bundle on each fiber, satisfies 
\[ \omega_\pi = \pi^\ast (\pi_\ast \omega_\pi) = \pi^\ast \mathbb{E} .\] 
In particular,
\[ c_1(\omega_\pi) = \pi^\ast \lambda .\]

Note that the choice of $\g{L}$ is canonical only up to a twist $\g{L} \mapsto \g{L}\otimes\pi^\ast \alpha$, for some $\alpha\in\Pic(\nik)$. A polarized Nikulin surface $(X,L)$ comes equipped with a distinguished set of $8$ smooth rational curves $R_1,\ldots,R_8$ which form, together with
\[ e = \frac{1}{2} (R_1+\cdots+R_8) ,\] 
a basis of the Nikulin lattice $\f{N}$. We consider the universal bundle
\[ \g{E} \in \Pic(\g{X}) \] 
associated to the distinguished class $e\in \Pic(X)$. Again, the choice of $\g{E}$ is canonical only up to a twist $\g{E} \mapsto \g{E} \otimes \pi^\ast \beta$, for some $\beta \in \Pic(\nik)$.

Following \cite[\S 4]{MOP17} we consider in the Chow ring of $\nik$ the kappa classes
\[ \kappa_{a,b,c} = \pi_\ast \left( c_1(\g{L})^a \cdot c_1(\g{E})^b \cdot c_2( \g{T}_\pi )^c \right) \in \mathrm{CH}^{a+b +2c -2} (\nik). \]
for some non-negative integers $a,b,c$. The four classes in codimension zero are readily computed by restriction to a fiber $(X,L)$ as follows:
\begin{align*}
\kappa_{2,0,0} &= L^2 = 2g-2 ,\\
\kappa_{0,2,0} &= e^2 = -4 , \\
\kappa_{1,1,0} &= L\cdot e = 0 ,\\
\kappa_{0,0,1} &= \chi_{\textit{top}}(X) = 24.
\end{align*}

In codimension one, we get six classes:
\[ \kappa_{3,0,0}, \, \, \kappa_{0,3,0}, \, \, \kappa_{2,1,0}, \, \, \kappa_{1,2,0}, \, \, \kappa_{1,0,1}, \, \, \kappa_{0,1,1} . \]
By twisting $(\g{L},\g{E}) \mapsto (\g{L}\otimes\pi^\ast \alpha , \g{E}\otimes \pi^\ast\beta)$, these classes will change in general. However, it is straightforward to verify that the following linear combinations are independent from these twists, and thus are canonically defined in $\mathrm{CH}^1(\nik)$,

\begin{align*}
\gamma_0 = \kappa_{3,0,0} - \frac{(g-1)}{4} \kappa_{1,0,1} \, , \quad \gamma_2 = \kappa_{1,0,1} + 6 \kappa_{1,2,0} \, , \\
\gamma_1 = \kappa_{2,1,0} - \frac{(g-1)}{12} \kappa_{0,1,1} \, , \quad \gamma_3 =\kappa_{0,1,1} + 2 \kappa_{0,3,0} \, .
\end{align*}

As in the previous sections, we write
\[ g = 2k^2 + p, \quad k\geq1, \quad 0\leq p < 4k+2 \]
and denote by $g_m = g-2m^2$ the genus of the line bundle $L_m = L\otimes e^{-m}$.

\bprop\label{GRR} Let $0\leq m\leq k$ and assume $p\geq 2$ if $m=k$. Following the above notations, let $\pi \colon \g{X} \to \nik$ be the universal polarized Nikulin surface of genus $g = 2k^2+p$,  with universal bundles $\g{L}$ and $\g{E}$, and denote by $\g{L}_m=\g{L}\otimes \g{E}^{\otimes -m}$. For any $n\geq1$, the sheaf $\g{U}_{n,m} := \pi_\ast ( \g{L}_m^{\otimes n})$ is locally free of rank $n^2(g_m -1)+2$ on $\nik$ with first Chern class
\begin{align*}
c_1(\g{U}_{n,m}) = & \, \, \, \frac{n^3}{6} [ \kappa_{3,0,0} -3m \kappa_{2,1,0} + 3m^2 \kappa_{1,2,0}  -m^3 \kappa_{0,3,0} ] +  \\
& + \frac{n}{12} [\kappa_{1,0,1} -m \kappa_{0,1,1} ] - [1+\frac{n^2}{2} (g_m-1) ] \lambda.
\end{align*}
\nprop

\begin{proof}
For $(X,L)\in \nik$, we have $h^0(X, L_m^{\otimes n}) = n^2(g_m -1)+2$ and $h^i(X,L_m^{\otimes n})=0 $ for $i=1,2$, whence $\pi_\ast \g{L}_m^{\otimes n}$ is a locally free sheaf of rank $n^2(g_m -1)+2$ and higher direct images vanish $R^i\pi_\ast \g{L}_m^{\otimes n}=0$. By the Grothendieck-Riemann-Roch formula,
\[ \ch \left(  \pi_\ast  \g{L}_m^{\otimes n} \right) = \pi_\ast \left\lbrace \ch (\g{L}_m^{\otimes n} ) \cdot \td(\g{T}_\pi) \right\rbrace . \]

Since $\pi$ is of relative dimension 2, the first Chern class of $\pi_\ast \g{L}_m^{\otimes n}$ equals the push-forward of the degree 3 part of the product on the right-hand side,
\[ c_1(\pi_\ast \g{L}_m^{\otimes n}) = \pi_\ast \left\lbrace \ch (\g{L}_m^{\otimes n} ) \cdot \td(\g{T}_\pi) \right\rbrace_{3} . \]
Recall that the Chern character and Todd class can be expressed in terms of Chern classes as follows:
\begin{align*}
\ch & = \rk + c_1 + \frac{1}{2} (c_1^2 -2c_2) + \frac{1}{6}(c_1^3 -3c_1c_2 +3c_3) + \cdots \\
\td & = 1 + \frac{c_1}{2} + \frac{1}{12}(c_1^2 + c_2) + \frac{1}{24}(c_1c_2) + \cdots
\end{align*}

Now, since $c_1(\omega_\pi) = \pi^\ast \lambda$, we get
\begin{align*}
\pi_\ast( c_1(\omega_\pi) \cdot c_2(\g{T}_\pi))  &= \lambda \cdot \pi_\ast c_2(\g{T}_\pi) = 24 \lambda, \\
\pi_\ast( c_1(\g{L}_m)^2 \cdot c_1(\omega_\pi)) &= \pi_\ast(c_1(\g{L}_m)^2)\cdot \lambda = (2g_m -2)\lambda , \\
\pi_\ast( c_1(\g{L}_m) \cdot c_1(\omega_\pi)^2) &= \underbrace{\pi_\ast( c_1(\g{L}_m))}_{=0} \cdot \lambda^2 =0 .
\end{align*}

Using these and $c_1(\g{T}_\pi)= -c_1(\omega_\pi)$, we get that $c_1(\pi_\ast \g{L}_m^{\otimes n})$  equals
\begin{align*}
 & \pi_\ast\biggl\lbrace \left( 1 + n c_1(\g{L}_m) + \frac{n^2}{2} c_1(\g{L}_m)^2 + \frac{n^3}{6} c_1(\g{L}_m)^3 \right) \cdot  \\
& \phantom{ == } \cdot \left( 1 + \frac{c_1(\g{T}_\pi)}{2} + \frac{c_1(\g{T}_\pi)^2 + c_2(\g{T}_\pi) }{12} + \frac{ c_1(\g{T}_\pi)c_2(\g{T}_\pi) }{24} \right) \biggl\rbrace_3 \\
&= \frac{n^3}{6} \pi_\ast(c_1(\g{L}_m)^3) - \frac{n^2}{4} (2g_m-2)\lambda + \frac{n}{12} \pi_\ast(c_1(\g{L}_m) c_2(\g{T}_\pi) )  - \lambda.  
\end{align*}
Substituting $c_1(\g{L}_m)= c_1(\g{L}) - mc_1(\g{E})$ and expanding the last expression yields to the claimed formula.
\end{proof}

We are now ready to prove Theorem \ref{THM-class-of-Dm} from the Introduction.

\begin{proof}[Proof of Theorem.\,\ref{THM-class-of-Dm} ]
We consider the morphism of vector bundles over $\nik$,
\[ \psi_m \colon \, S^2 \g{U}_{1,m} \to \g{U}_{2,m} . \]
The divisor $\f{D}^{\rk (4)}_m$ is by definition the locus where $\ker(\psi_m)$ contains quadrics of rank at most four. Applying the formula from \cite[Thm.\,1.1]{FR}, we get
\[ \left[ \f{D}^{\rk (4)}_m \right] =  A^{g_m-3}_{g_m+1} \left( c_1( \g{U}_{2,m} ) - 2 \left( \frac{ \rk \g{U}_{2,m} }{ \rk \g{U}_{1,m} } \right) c_1( \g{U}_{1,m}) \right) , \]
where the coefficient corresponds to the number $A^r_e$,
\[
A^r_e = \frac{ \binom{e}{r} \binom{e+1}{r-1} \cdots \binom{e+r-1}{1} }{ \binom{1}{0} \binom{3}{1} \binom{5}{2} \cdots \binom{2r-1}{r-1} }
\]
which represents the degree of the variety of symmetric $e\times e$ matrices of corank  at least $r$ inside the projective space of symmetric $e\times e$ matrices.

By Proposition \ref{GRR}, the class
\[ c_1( \g{U}_{2,m} ) - 2 \left( \frac{ \rk \g{U}_{2,m} }{ \rk \g{U}_{1,m} } \right) c_1( \g{U}_{1,m})  \] 
is equal to
\begin{align*}
& \frac{2}{ g_m + 1 } [ \kappa_{3,0,0} -3m \kappa_{2,1,0} + 3m^2 \kappa_{1,2,0}  -m^3 \kappa_{0,3,0} ]  \\ 
& - \frac{ g_m-1}{ 2(g_m+1) } [\kappa_{1,0,1} -m \kappa_{0,1,1} ]  + (2g_m-1) \lambda.
\end{align*}

Plugging into the latter expression the following straightforward identities:
\[
\begin{cases}
2\gamma_0 + m^2\gamma_2  = 2\kappa_{3,0,0} - \frac{g_m-1}{2} \kappa_{1,0,1} +6m^2 \kappa_{1,2,0} , \\
-m^3 \gamma_3 -6m \gamma_1  = -2m^3\kappa_{0,3,0} + \frac{g_m-1}{2} m \kappa_{0,1,1} -6m\kappa_{2,1,0} ,
\end{cases} 
\]
one gets

\begin{displaymath}
\frac{ 2 }{ g_m +1 } \gamma_0 - \frac{ 6m }{ g_m +1 } \gamma_1 + \frac{ m^2 }{ g_m +1 } \gamma_2 - \frac{ m^3 }{ g_m +1 } \gamma_3 + (2g_m-1) \lambda . 
\end{displaymath} \phantom{.}

The result follows.
\end{proof}

\bibliography{bigbib}
\bibliographystyle{plain}

\end{document}